# Estimation of AR and ARMA models by stochastic complexity

**Ciprian Doru Giurcăneanu**[1,*] **and Jorma Rissanen**[2,†,‡]

*Tampere University of Technology, and Technical University of Tampere and Helsinki, and Helsinki Institute for Information Technology*

**Abstract:** In this paper the stochastic complexity criterion is applied to estimation of the order in AR and ARMA models. The power of the criterion for short strings is illustrated by simulations. It requires an integral of the square root of Fisher information, which is done by Monte Carlo technique. The stochastic complexity, which is the negative logarithm of the Normalized Maximum Likelihood universal density function, is given. Also, exact asymptotic formulas for the Fisher information matrix are derived.

## 1. Introduction

The negative logarithm of the *NML* (Normalized Maximum Likelihood) universal model, called the *stochastic complexity*, provides a powerful criterion for estimation of the model structure such as the optimal collection of the regressor variables in the linear quadratic regression problem, [19], especially for small amounts of data. It involves the integral of the square root of the Fisher information, which is easy to calculate when the regressor matrix does not depend on the parameters. While modeling gaussian time series with AR models are instances of linear quadratic regression problems their order estimation poses trouble with the stochastic complexity for the reason that the regressor matrix is determined by the parameters, and the Fisher information is not constant. The same problem of course is also with the ARMA models, which have the additional difficulty of calculation of the maximum likelihood parameters.

In this paper we resort to Monte Carlo integration to overcome the problem posed by the nonconstant Fisher information and study by simulations the efficiency of the resulting order estimation criterion. Although exact formulas exist for the Fisher information matrix they are quite cumbersome to evaluate, and we consider asymptotic simplifications. This may run against the intent of getting a criterion for small amounts of data, but the asymptotic estimates appear to be good enough, and the resulting criterion for the short data sequences created is still superior among the competing criteria such as the *BIC* [20], which is equivalent with a

[1]Institute of Signal Processing, Tampere University of Technology, P.O. Box 553, FIN-33101 Tampere, Finland, e-mail: `ciprian.giurcaneanu@tut.fi`

[2]140 Teresita Way, Los Gatos, CA 95032, USA, e-mail: `jrrissanen@yahoo.com`

*The work of C.D. Giurcăneanu was supported by the Academy of Finland, project No. 213462 (Finnish Centre of Excellence program (2006–2011).

†J. Rissanen is affiliated with the Technical Universities of Tampere and Helsinki, and Helsinki Institute for Information Technology, Finland.

‡Corresponding author.

crude asymptotic version of the *MDL* criterion [15], and a recently suggested one, *KICC* [21], or bias corrected Kullback-Leibler criterion.

We describe below the *NML* model for AR and ARMA class of models, and discuss its optimality properties. We also derive in the Appendix the asymptotic form of the Fisher information matrix for the general ARMA class of models.

## 2. Normalized maximum likelihood model

We consider the ARMA model:

$$(1) \quad y_t + \sum_{i=1}^{n} a_i y_{t-i} = e_t + \sum_{j=1}^{m} b_j e_{t-j},$$

where $e_t$ is zero-mean white Gaussian noise of variance $\sigma^2$. The integers $m, n$ are nonnegative, and all coefficients $a_i$ and $b_j$ are real-valued. We can equivalently write $y_t = \frac{B(q)}{A(q)} e_t$, where $B(q) = 1 + b_1 q^{-1} + \cdots + b_m q^{-m}$, $A(q) = 1 + a_1 q^{-1} + \cdots + a_n q^{-n}$, and $q^{-1}$ is the unit delay operator. We will use the notation ARMA(n,m) for the class of the normal density functions $\{f(y^N; \theta)\}$ defined by such processes, where $\theta = (a_1, \ldots, a_n, b_1, \ldots, b_m, \sigma^2)$, the parameters ranging over a subset of $\Re^k$, where $k = n + m + 1$. Let $\hat{\theta}(y^N)$ denote the maximum likelihood estimates of the parameters $\theta$.

In order to define the range of the parameters properly we need to consider another equivalent parametrization in terms of the roots of the two polynomials

$$(2) \quad \prod_{i=1}^{n}(1 - g_i q^{-1}) y_t = \prod_{j=1}^{m}(1 - h_j q^{-1}) e_t,$$

together with the noise variance $\sigma^2$. We denote by $g_i$ the zeros of $A(q)$ and by $h_j$ the zeros of $B(q)$. There are no repeated poles or zeros nor pole-zero cancellations. We specify in the Appendix exactly the further restrictions on the type of the zeros but for now let the same symbol $\theta$ denote the new parameters ranging over $\Theta \subset \Re^k$.

Consider the *NML* density function, [3],[18],

$$\hat{f}(y^N; n, m) = \frac{f(y^N; \hat{\theta}(y^N))}{C_{k,n}},$$

where

$$C_{k,n} = \int_{x^N: \hat{\theta}(x^N) \in \Omega} f(x^N; \hat{\theta}(x^N)) dx^N$$
$$= \int_{\hat{\theta} \in \Omega} g(\hat{\theta}; \hat{\theta}) d\hat{\theta},$$

and $g(\hat{\theta}; \theta)$ denotes the density function on the statistic $\hat{\theta}$ induced by $f(y^N; \theta)$. In the equation above, we use the identity $f(x^N; \hat{\theta}(x^N), \theta) = f(x^N | \hat{\theta}(x^N); \theta) g(\hat{\theta}(x^N); \theta)$, that is integrated first over $x^N$ at the point $\hat{\theta}(x^N) = \hat{\theta} = \theta$ kept fixed, which gives unity, and then over $\hat{\theta}$.

Under the main assumption that the convergence in distribution by the Central Limit Theorem applies to the ML estimates, the stochastic complexity, $L(y^N;$

$n, m) = \ln 1/\hat{f}(y^N; n, m)$, is given by

$$(3) \quad L(y^N; n, m) = -\ln f(y^N; \hat{\theta}(y^N)) + \frac{k}{2}\ln\frac{N}{2\pi} + \ln\int_\Theta |\mathbf{J}(\theta)|^{1/2}d\theta + o(1),$$

where $\Theta$ denotes the parameter space, and $\mathbf{J}(\theta)$ is the Fisher information matrix [18]. The rate of convergence o(1) is determined by the convergence of the ML estimates to the normal density function.

To get a criterion for the structure in general we ought to add the code length needed to encode the structure, but here we take the simple case where the structure consists of a few first coefficients of the ARMA model, whose code length is much shorter than the stochastic complexity and ignored. (If $k$ is not small, we can use the estimate $L(k) = \ln k + 2\ln\ln k$.)

The *NML* model has the following two optimality properties, which justify its name:

(1) It is the unique solution $\hat{f} = \hat{g} = \hat{q}$ to the following maxmin problem

$$\max_g \min_q E_g \log \frac{f(y^N; \hat{\theta}(y^N))}{q(y^N)},$$

where $g$ and $q$ range over any sets that include $\hat{f}$. Notice that the logarithm of the ratio is the difference between the ideal code length $\log 1/q$ and the unattainable lower bound for any code length in the ARMA class.

(2) If the data generating distribution $g$ is restricted to the ARMA class, the mean of the stochastic complexity with respect to the model $\theta$ cannot be beaten by any model what so ever, except for $\theta$ in a set whose volume goes to zero as $N$ grows.

## 3. Linear regression with constant regressor matrix

Before discussing the AR models we illustrate the stochastic complexity criterion for linear quadratic regression with constant Fisher information by comparing it with the *BIC* and the *KICC* criteria in a simple polynomial fitting problem for small amounts of data.

For linear regression with a constant regressor matrix $\mathbf{X} = \{x_{it}\}$ the stochastic complexity criterion takes the form, [19],

$$\min_{\gamma \in \Gamma}\{(N-k)\ln\hat{\tau} + k\ln\hat{R} + (N-k-1)\ln\frac{1}{n-k} - (k-1)\ln k\}.$$

The index $\gamma = i_1, \ldots, i_k$, consists of the indices of the rows $\bar{\mathbf{x}}_\mathbf{i}$ of the $k \times n$ regressor matrix included in the linear combination

$$y_t = \sum_{i \in \gamma} \beta_i \bar{x}_{it} + e_t, \ t = 1, \ldots, N,$$

$\hat{\tau}$ is the minimized squared error per symbol, and $\hat{R} = \frac{1}{n}\hat{\boldsymbol{\beta}}^\top \mathbf{X}_\gamma \mathbf{X}_\gamma^\top \hat{\boldsymbol{\beta}}$, where $\mathbf{X}_\gamma$ is the $k \times n$ submatrix of $\mathbf{X}$ consisting of the retained rows.

Notice that there are no hyper parameters defining the range of the parameters $\beta_i$ and $\tau$. They have been renormalized away.

TABLE 1
*Order estimation of the polynomial model in Example 1. The true order is $k = 3$. For each criterion, the probability of correct estimation of the order is computed from $10^5$ runs. Also shown is the probability of overestimation of the polynomial order ($4 \leq \hat{k} \leq 10$). The probability of underestimation ($0 \leq \hat{k} \leq 2$) is almost zero for all analyzed criteria. The best result for each sample size $N$ is represented with bold font.*

| Order | Criterion | Sample size($N$) | | | | | | | | |
|---|---|---|---|---|---|---|---|---|---|---|
| | | 25 | 30 | 40 | 50 | 60 | 70 | 80 | 90 | 100 |
| **$\hat{k} = k$** | NML | **0.94** | **0.95** | **0.96** | **0.97** | **0.97** | **0.97** | **0.98** | **0.98** | **0.98** |
| | BIC | 0.79 | 0.84 | 0.89 | 0.91 | 0.93 | 0.94 | 0.95 | 0.95 | 0.95 |
| | KICC | 0.93 | 0.92 | 0.91 | 0.91 | 0.90 | 0.90 | 0.90 | 0.90 | 0.89 |
| $\hat{k} > k$ | NML | 0.06 | 0.05 | 0.04 | 0.03 | 0.03 | 0.03 | 0.02 | 0.02 | 0.02 |
| | BIC | 0.21 | 0.16 | 0.11 | 0.09 | 0.07 | 0.06 | 0.05 | 0.05 | 0.05 |
| | KICC | 0.07 | 0.08 | 0.09 | 0.09 | 0.10 | 0.10 | 0.10 | 0.10 | 0.11 |

**Example 1.** We discuss an example of polynomial fitting considered in [21] to investigate the performances of a model selection criterion called *KICC*. It is obtained by an application of a bias correction to KIC (Kullback Information Criterion), [6], and it is recommended to be used in linear regression problems when the sample size is small. The underlying signal is generated by a third-order polynomial model $\tilde{y} = x^3 - 0.5x^2 - 5x - 1.5$, where the points $x_1, \ldots, x_N$ are chosen to be uniformly distributed in $[-3, 3]$. The measurements $y_1, \ldots, y_N$ are obtained by addition to $\tilde{y}_i$ zero-mean white Gaussian noise, whose variance is selected such that the signal-to-noise ratio is SNR=10 dB. For each number of data points $N$, between 25 and 100, $10^5$ different realizations are produced, to which polynomials of degree $0, 1, \ldots, 10$ are fitted with the least squares method.

The estimates of the order of the polynomial obtained with the *NML*, *BIC* and *KICC* criteria are in Table 1. We have restricted our investigations only to these three criteria, because in [21] *KICC* was shown to outperform other six estimation criteria for $N = 25$ and $N = 30$. We see in the table that *NML* criterion performs better than *BIC* and *KICC* in all the cases studied. Observe that the number of correct estimations produced by *KICC* generally declines when more measurements are available, while the *BIC* and the *NML* results improve with increasing $N$. For example, *KICC* compares favorable with *BIC* for $N = 25$, but the situation is reversed for $N = 100$.

## 4. AR models

The likelihood density function for an AR model is given by

$$f(y^N; \theta) = \frac{1}{(2\pi\sigma^2)^{N/2}} e^{-\frac{1}{2\sigma^2} \sum_{t=1}^{N}(y_t + a_1 y_{t-1} + \cdots + a_n y_{t-n})^2},$$

where we put $y_t = 0$ for $t < 1$. The maximized likelihood is $\frac{1}{(2\pi e \hat{\sigma}^2)^{N/2}}$, where $\hat{\sigma}^2$ is the minimized sum per symbol $\hat{\sigma}^2 = \frac{1}{N} \sum_{t=1}^{N}(y_t + \hat{a}_1 y_{t-1} + \cdots + \hat{a}_n y_{t-n})^2$. The *NML* criterion (3) has now the expression

(4) $\quad L(y^N; n) = \frac{N}{2} \ln(2\pi e \hat{\sigma}^2) + \frac{n+1}{2} \ln \frac{N}{2\pi} + \ln \int_\Theta |\mathbf{J}(\theta)|^{1/2} d\theta + o(1).$

The Fisher information matrix is given by $\begin{bmatrix} \mathbf{R}_{zz} & 0 \\ 0 & 1/(2\sigma^4) \end{bmatrix}$, where

$$\mathbf{R}_{zz} = \begin{bmatrix} r_0 & r_1 & \cdots & r_{n-1} \\ r_1 & r_0 & \cdots & r_{n-2} \\ \vdots & \vdots & \ddots & \vdots \\ r_{n-1} & r_{n-2} & \cdots & r_0 \end{bmatrix},$$

and $r_i = E[z_t z_{t-i}]$ denote the covariances of the process $z_t = y_t/\sigma$ [9, 10]. Applying the formula in [12] for the parameters transformation and the well-known Vieta's formulae, it is easy to calculate the Fisher information matrix for the parameter set given by the model poles $g = (g_1, g_2, \ldots, g_n)$ and the noise variance $\sigma^2$.

Remark in (4) that the integral term makes the most important difference between the expression for the stochastic complexity and the BIC criterion. The integral has a lot of structural information which BIC lacks, and it generally increases with $n$, because the determinant increases.

We note that the contribution of the $\sigma^2$ to the integral is decoupled by the contribution of the other parameters. Consequently we ignore for all the AR models the contribution of $\sigma^2$ because we do not have any "natural" finite limits for the range of $\sigma^2$. The constrain to have a stable model restricts the domain of the magnitudes of the poles to be a hypercube.

Apart from the AR(1) case for which the integral in (3) can be found in a closed form, $\int_{-1}^{1} \frac{1}{\sqrt{1-g^2}} dg = \pi$, the evaluation of the integral will be done by the Monte Carlo technique. To be more precise we use Sobol' sequences [14] to perform the Monte Carlo integration for AR(n) models with $1 \leq n \leq 6$. For these values all poles are complex if $n$ is even, and exactly one pole is real-valued if $n$ is odd, which can be taken advantage of in calculating the form of the information matrix.

Our Matlab implementation is based on the algorithm described on p. 312 in [14] and the code publicly available at [1]. We perform the Monte Carlo integration for various AR models with $M$ integration points. But first, to test the accuracy we use the known result for the AR(1) model. Table 2 shows the fractional error obtained when $M = 10^5$ and $M = 10^6$. For models with larger order, we report the value $\Delta = |\hat{I}_{10^7} - \hat{I}_{10^6}|/\hat{I}_{10^7}$, where $\hat{I}_M$ denotes the Monte Carlo evaluation of $\int_\Theta |\mathbf{J}(g)|^{1/2} d\theta$ calculated from $M$ integration points. We show in Table 2 the results on $\Delta$ since it is known for Monte Carlo integration with Sobol' sequences that the fractional error decreases with the number of samples as $(\ln M)^n/M$ [14].

TABLE 2
*Monte Carlo results for the integral term in the stochastic complexity formula (4) for autoregressive models. For the AR(1) model the fractional error is reported.*

| $M$ | $\hat{I}_M$ | Fractional error or $\Delta$ |
|---|---|---|
| AR(1) | | |
| $10^5$ | 3.131956 | 0.003067 |
| $10^6$ | 3.138952 | 0.000840 |
| AR(2) - pure complex poles | | |
| $10^6$ | 42.06 | - |
| $10^7$ | 47.41 | 0.11 |
| AR(3) - one real-valued pole | | |
| $10^6$ | 122.67 | - |
| $10^7$ | 137.73 | 0.11 |
| AR(4) - pure complex poles | | |
| $10^6$ | 1069.66 | - |
| $10^7$ | 1358.84 | 0.21 |
| AR(5) - one real-valued pole | | |
| $10^6$ | 3733.59 | - |
| $10^7$ | 8307.55 | 0.55 |
| AR(6) - pure complex poles | | |
| $10^6$ | 23164.39 | - |
| $10^7$ | 35981.48 | 0.36 |

**Example 2.** We evaluate the capabilities of *NML*, *BIC* and *KICC* criteria for estimating the order of AR models. The *NML* criterion is calculated with formula (4), where the value of the integral term for $n > 1$ is the one from Table 2 computed with $M = 10^7$ integration points. We extend our experimental framework by considering another information theoretic criterion, namely the predictive least squares criterion *PLS*, [16].

Figure 1 outlines the simulation procedure used in Example 2, and the estimation results are shown in Tables 3-4.

Note that the evaluation of the various criteria for order estimation requires the estimate of noise variance for each order between one and six. Moreover, for the *PLS* criterion the computation of the prediction errors must be performed for each order and for each sample point. To reduce the computational burden, we resort to the fast implementation of the prewindowed estimation method based on predictive lattice filters [8], [22].

Observe in Table 3 that the *NML* criterion compares favorably with all the other criteria when the sample size is at least 50. For the smallest amount of data

---

<u>For</u> the model order $n \in \{1, 2, 3\}$,
    For each order estimation criterion $\mathcal{C}$ and for each sample size $N$,
        $N \in \{25, 50, 100, 200\}$, initialize with zero two counters:
        $\mathcal{N}^c_{N,\mathcal{C}}$ for correct estimations and $\mathcal{N}^o_{N,\mathcal{C}}$ for over-estimations.
    <u>Repeat</u> the following steps 1000 times:
        Generate independently the entries of $\mathcal{P}_\mu$ as outcomes of $\mathcal{U}\,[(0.8, 1)]$,
        and the entries of $\mathcal{P}_\phi$ as outcomes of $\mathcal{U}\,[(0, \pi)]$.
        If $n$ is odd, generate the unique entry of $\mathcal{P}_\rho$
        according to $\mathcal{U}\left[(-1, -0.8) \bigcup (0.8, 1)\right]$.
    <u>Repeat</u> the following steps 1000 times:
        Simulate a time series with 300 entries for the AR(n) process
            whose poles are given by $\mathcal{P}_\mu$, $\mathcal{P}_\phi$, $\mathcal{P}_\rho$.
            Use null initial conditions and $\sigma^2 = 1$.
            Discard the first 100 entries of the time series and
            dub $z$ the vector formed with the rest of 200 measurements.
        <u>For</u> each sample size $N \in \{25, 50, 100, 200\}$,
            Choose $y^N = [z_1, \ldots, z_N]^\top$. Apply each criterion $\mathcal{C}$
            to estimate the model order $\hat{n}_{N,\mathcal{C}}$ from $y^N$ data,
            under the hypothesis $\hat{n}_{N,\mathcal{C}} \in \{1, \ldots, 6\}$.
            If $\hat{n}_{N,\mathcal{C}} = n$, then increment $\mathcal{N}^c_{N,\mathcal{C}}$.
            If $\hat{n}_{N,\mathcal{C}} > n$, then increment $\mathcal{N}^o_{N,\mathcal{C}}$.
        End
    End
    End
    Calculate the probability of correct estimation $\hat{p}^c_{N,\mathcal{C}} = \mathcal{N}^c_{N,\mathcal{C}}/10^6$,
        and the probability of over-estimation $\hat{p}^o_{N,\mathcal{C}} = \mathcal{N}^o_{N,\mathcal{C}}/10^6$ for the model order.
End

---

FIG 1. *The simulation procedure applied in Example 2. The notation $\mathcal{U}[\cdot]$ is used for the uniform distribution.*

the asymptotic calculation of the Fisher information does not seem to be accurate enough. In most of the cases *BIC* is ranked the second after the *NML*, and the results of *KICC* do not improve when the sample size $N$ is increased. For all criteria the performances decline for the larger values of the model order, which is clear because there is more to learn. Notice the moderate performances of the *PLS* criterion. We mention that another comparative study [7] also reports the moderate capabilities of *PLS* on estimating the order of AR models. This is to be expected since the *PLS* criterion is based on the estimates of the parameters which are shaky for small amounts of data.

## 5. ARMA models

The density function for ARMA models, (1), depends on how the initial values of $y$ are related to the inputs $e$. A simple formula results if we put $y_i = e_i = 0$ for $i \leq 0$. Then the linear spaces spanned by $y^t$ and $e^t$ are the same. Let $\hat{y}_{t+1|t}$ be the orthogonal projection of $y_{t+1}$ on the space spanned by $y^t$. We have the recursion

$$(5) \qquad \hat{y}_{t+1|t} = \sum_{i=1}^{m} b_i(y_{t-i+1} - \hat{y}_{t-i+1|t-i}) - \sum_{i=1}^{n} a_i y_{t-i+1},$$

where $\hat{y}_{1|0} = 0$. With more general initial conditions the coefficients $b_i$ in (5) will depend on t; see for instance [17]. The likelihood function of the model is then

$$(6) \qquad f(y^N; \theta, \sigma^2) = \frac{1}{(2\pi\sigma^2)^{N/2}} e^{-\frac{1}{2\sigma^2} \sum_{t=1}^{N}(y_t - \hat{y}_{t|t-1})^2}.$$

The maximized likelihood is $\frac{1}{(2\pi e\hat{\sigma}^2)^{N/2}}$, where $\hat{\sigma}^2 = \min_{a_1,\ldots,a_n,b_1,\ldots,b_m} \frac{1}{N}\sum_{t=1}^{N}(y_t - \hat{y}_{t|t-1})^2$. The *NML* criterion (3) is then given by

$$(7) \qquad L(y^N; n, m) = \frac{N}{2}\ln(2\pi e\hat{\sigma}^2) + \frac{n+m+1}{2}\ln\frac{N}{2\pi} + \ln\int_\Theta |\mathbf{J}(\theta)|^{1/2}d\theta + o(1),$$

In Appendix we elaborate on the computation of the integral term for the *NML* criterion, and the results are applied to the selection for ARMA models in the following example.

TABLE 3
*Example 2 - the probability of correct estimation of the AR order. The best result for each sample size N is represented with bold font.*

| AR model order | Criterion | Sample size ($N$) | | | |
|:---:|:---:|:---:|:---:|:---:|:---:|
| | | 25 | 50 | 100 | 200 |
| $n=1$ | NML | **0.99** | **0.99** | **1.00** | **1.00** |
| | BIC | 0.93 | 0.95 | 0.97 | 0.98 |
| | KICC | 0.95 | 0.93 | 0.91 | 0.90 |
| | PLS | 0.89 | 0.92 | 0.95 | 0.97 |
| $n=2$ | NML | 0.72 | **0.85** | **0.87** | **0.88** |
| | BIC | 0.79 | **0.85** | **0.87** | 0.87 |
| | KICC | **0.82** | 0.83 | 0.80 | 0.78 |
| | PLS | 0.49 | 0.59 | 0.66 | 0.71 |
| | NML | 0.49 | **0.74** | **0.83** | **0.84** |
| $n=3$ | BIC | **0.52** | 0.71 | 0.78 | 0.79 |
| | KICC | 0.51 | 0.71 | 0.73 | 0.69 |
| | PLS | 0.26 | 0.39 | 0.47 | 0.53 |

TABLE 4
*Example 2 - the probability to over-estimation of the order of AR models. The smallest overestimation probability for each sample size N is represented with bold font.*

| AR model order | Criterion | Sample size (N) | | | |
|---|---|---|---|---|---|
| | | 25 | 50 | 100 | 200 |
| $n = 1$ | NML | **0.01** | **0.01** | **0.00** | **0.00** |
| | BIC | 0.07 | 0.05 | 0.03 | 0.02 |
| | KICC | 0.05 | 0.07 | 0.09 | 0.10 |
| | PLS | 0.11 | 0.08 | 0.05 | 0.03 |
| $n = 2$ | NML | 0.07 | **0.09** | **0.11** | **0.12** |
| | BIC | 0.10 | 0.11 | 0.12 | 0.13 |
| | KICC | **0.06** | 0.14 | 0.20 | 0.22 |
| | PLS | 0.20 | 0.19 | 0.17 | 0.15 |
| $n = 3$ | NML | **0.01** | **0.03** | **0.06** | **0.12** |
| | BIC | 0.07 | 0.09 | 0.12 | 0.18 |
| | KICC | 0.03 | 0.10 | 0.20 | 0.29 |
| | PLS | 0.21 | 0.22 | 0.23 | 0.23 |

TABLE 5
*Results of model selection for the ARMA models in Example 3. The counts indicate for 1000 runs the number of times the structure of the model was correctly estimated by each criterion, from the set $\{ARMA(n,m) : n, m \geq 1, n + m \leq 6\}$. The best result for each sample size N is represented with bold font.*

| ARMA model | Criterion | Sample size (N) | | | | |
|---|---|---|---|---|---|---|
| | | 25 | 50 | 100 | 200 | 400 |
| $n = 1, m = 1$ | NML | 700 | **812** | **917** | **962** | **989** |
| $a_1 = -0.5$ | BIC | 638 | 776 | 894 | 957 | 983 |
| $b_1 = 0.8$ | KICC | **717** | 740 | 758 | 745 | 756 |
| $n = 2, m = 1$ | NML | **626** | **821** | **960** | **991** | **994** |
| $a_1 = 0.64, a_2 = 0.7$ | BIC | 532 | 740 | 898 | 961 | 978 |
| $b_1 = 0.8$ | KICC | 586 | 727 | 810 | 846 | 849 |
| $n = 1, m = 1$ | NML | 851 | **887** | **918** | **931** | **961** |
| $a_1 = 0.3$ | BIC | 766 | 804 | 856 | 903 | 942 |
| $b_1 = 0.5$ | KICC | **860** | 764 | 654 | 614 | 577 |

**Example 3.** We calculate the structure of ARMA models for data generated by three different processes, which also were used in [11]. For each model, the true structure and the coefficients are given in Table 5, where we show the estimation results for 1000 runs. In all experiments we have chosen the variance of the zero-mean white Gaussian noise to be $\sigma^2 = 1$. We mention that, similarly with the experiments on the autoregressive models each data set $y^N$ was obtained after discarding the first 100 generated measurements. This is to eliminate the effect of the initial conditions. There exist different methods for estimation of ARMA models. We selected the one implemented in Matlab as armax function by Ljung, which is well described in his book [13].

**Appendix: The asymptotic Fisher information matrix**

We focus on the computation of the integral term in equation (7). The model is assumed to be stable and minimum phase, which means that in (2) the roots for both $B(q)$ and $A(q)$ are inside the open unit disc. Assume that $n_1$ zeros of $A(q)$ and $m_1$ zeros of $B(q)$ are real-valued. Then we have the inequalities $0 \leq n_1 \leq n$ and $0 \leq m_1 \leq m$. Because all coefficients of $A(q)$ and $B(q)$ are real-valued, the pure complex poles and zeros occur in complex conjugate pairs, and consequently the differences $n - n_1$ and $m - m_1$ are both even integers. For the pure complex poles and zeros we apply the parametrization in [5]:

$$g_{\ell+1} = g_\ell^* = |g_\ell|\exp(-i\phi_{g_\ell}), \quad \phi_{g_\ell} \in (0,\pi),\ \ell \in \{n_1+1, n_1+3, \ldots, n-1\},$$
$$h_{\ell+1} = h_\ell^* = |h_\ell|\exp(-i\phi_{h_\ell}), \quad \phi_{h_\ell} \in (0,\pi),\ \ell \in \{m_1+1, m_1+3, \ldots, m-1\},$$

where the symbol $*$ denotes the complex conjugate. The entries of the parameter vector $\theta$ are given by:

$$\begin{aligned}\theta = (&g_1, \ldots, g_{n_1}, \\ &|g_{n_1+1}|, \phi_{g_{n_1+1}}, \ldots, |g_{n-1}|, \phi_{g_{n-1}}, \\ &h_1, \ldots, h_{m_1}, \\ &|h_{m_1+1}|, \phi_{h_{m_1+1}}, \ldots, |h_{m-1}|, \phi_{h_{m-1}}, \\ &\sigma^2).\end{aligned}$$

For the sake of clarity we define the subsets of indices for the $\theta$ parameters:

$$\begin{aligned}\mathcal{P}_\rho &= \{1, 2, \ldots, n_1\} \\ \mathcal{P}_\mu &= \{n_1+1, n_1+3, \ldots, n-1\} \\ \mathcal{P}_\phi &= \{n_1+2, n_1+4, \ldots n\} \\ \mathcal{P} &= \mathcal{P}_\rho \bigcup \mathcal{P}_\mu \bigcup \mathcal{P}_\phi \\ \mathcal{Z}_\rho &= \{n+1, n+2, \ldots n+m_1\} \\ \mathcal{Z}_\mu &= \{n+m_1+1, n+m_1+3, \ldots, n+m-1\} \\ \mathcal{Z}_\phi &= \{n+m_1+2, n+m_1+4, \ldots, n+m\} \\ \mathcal{Z} &= \mathcal{Z}_\rho \bigcup \mathcal{Z}_\mu \bigcup \mathcal{Z}_\phi\end{aligned}$$

Based on (6) we use the following asymptotic expression for the log-likelihood function of the observations $y_1, \ldots, y_N$, [2], [9]:

$$\mathcal{L} = -\frac{1}{2\sigma^2}\sum_{t=1}^N e_t^2 - \frac{N}{2}\ln\sigma^2 + \text{constant}.$$

For all $u, v \in \{1, \ldots, m+n+1\}$, the $(u, v)$ entry of the Fisher information matrix is given by the formula [18]: $J_{u,v} = -\lim_{N\to\infty}\frac{1}{N}E[\frac{\partial^2\mathcal{L}}{\partial\theta_u\partial\theta_v}]$. Applying the results in [2] and [9], we obtain in a straightforward manner:

$$\begin{aligned}J_{n+m+1,n+m+1} &= 1/(2\sigma^4), \\ J_{u,n+m+1} = J_{n+m+1,v} &= 0\ \forall u, v \in \{1, \ldots, n+m\}.\end{aligned}$$

For the following calculations we use the identity $J_{u,v} = \lim_{N\to\infty}\frac{1}{N}E[\frac{\partial\mathcal{L}}{\partial\theta_u}\frac{\partial\mathcal{L}}{\partial\theta_v}]$. Consider first the case $u, v \in \mathcal{P}_\rho$. Simple calculations lead to

$$\frac{\partial e_t}{\partial\theta_u} = -\frac{q^{-1}}{1-\theta_u q^{-1}}e_t = -\sum_{p=1}^\infty \theta_u^{p-1}q^{-p}e_t,$$

and we obtain readily:

$$J_{u,v} = \frac{1}{N\sigma^4} E\left[\left(\sum_{t=1}^{N} e_t \sum_{p=1}^{\infty} \theta_u^{p-1} e_{t-p}\right)\left(\sum_{s=1}^{N} e_s \sum_{r=1}^{\infty} \theta_v^{r-1} e_{s-r}\right)\right]$$

(8)
$$= \frac{1}{N\sigma^4} \sum_{t=1}^{N} \sum_{p=1}^{\infty} (\theta_u \theta_v)^{p-1} E\left[e_t^2 e_{t-p}^2\right]$$

$$= \frac{1}{1 - \theta_u \theta_v}.$$

We conclude for $u, v \in \mathcal{P}_\rho \bigcup \mathcal{Z}_\rho$ that $J_{u,v} = \frac{\mathcal{S}_u \mathcal{S}_v}{1 - \theta_u \theta_v}$, where

$$\mathcal{S}_u = \begin{cases} -1, \, u \in \mathcal{P} \\ 1, \quad u \in \mathcal{Z} \end{cases}$$

Formula (8) was deduced in [4] for the case when all the poles and the zeros of the ARMA(n,m) model are real-valued. We evaluate next the entry $(u, v)$ of the Fisher information matrix for $u \in \mathcal{P}_\rho \bigcup \mathcal{Z}_\rho$ and $v \in \mathcal{P}_\mu \bigcup \mathcal{P}_\phi \bigcup \mathcal{Z}_\mu \bigcup \mathcal{Z}_\phi$. It is not difficult to prove that

$$\frac{\partial e_s}{\partial \theta_v} = \sum_{r=1}^{\infty} d_{v,r} e_{s-r} \,\, \forall s \in \{1, \ldots, N\},$$

where the coefficients $d_{v,r}$ are real-valued, [5]. Therefore

$$J_{u,v} = \frac{\mathcal{S}_u}{N\sigma^4} E\left[\left(\sum_{t=1}^{N} e_t \sum_{p=1}^{\infty} \theta_u^{p-1} e_{t-p}\right)\left(\sum_{s=1}^{N} e_s \sum_{r=1}^{\infty} d_{v,r} e_{s-r}\right)\right]$$

$$= \frac{\mathcal{S}_u}{N\sigma^4} \sum_{t=1}^{N} \sum_{p=1}^{\infty} \theta_u^{p-1} d_{v,p} E\left[e_t^2 e_{t-p}^2\right]$$

$$= \mathcal{S}_u \sum_{p=1}^{\infty} \theta_u^{p-1} d_{v,p}.$$

The following closed form expressions of $d_{v,p}$ are given in [5] for $v \in \mathcal{P}_\mu \bigcup \mathcal{Z}_\mu$:

$$d_{v,p} = \begin{cases} 2\mathcal{S}_v \cos\theta_{v+1}, & p = 1 \\ 2\mathcal{S}_v \frac{\theta_v^p \sin(p\theta_{v+1})\cos\theta_{v+1} - \theta_v^{p-1}\sin((p-1)\theta_{v+1})\theta_v}{\theta_v \sin\theta_{v+1}}, & p \geq 2 \end{cases}$$

The equations above lead to

$$J_{u,v} = 2\frac{\mathcal{S}_u \mathcal{S}_v}{\theta_u \theta_v} \frac{\cos\theta_{v+1}}{\sin\theta_{v+1}} \sum_{p=1}^{\infty} (\theta_u \theta_v)^p \sin(p\theta_{v+1})$$

$$-2\frac{\mathcal{S}_u \mathcal{S}_v}{\sin\theta_{v+1}} \sum_{p=1}^{\infty} (\theta_u \theta_v)^p \sin(p\theta_{v+1})$$

$$= 2\mathcal{S}_u \mathcal{S}_v \frac{\cos\theta_{v+1} - \theta_u \theta_v}{1 - 2\theta_u \theta_v \cos\theta_{v+1} + \theta_u^2 \theta_v^2},$$

for $u \in \mathcal{P}_\rho$ and $v \in \mathcal{P}_\mu \bigcup \mathcal{Z}_\mu$. Similarly for $v \in \mathcal{P}_\phi \bigcup \mathcal{Z}_\phi$ and $p \geq 1$, we have, [5],

$$d_{v,p} = -2\mathcal{S}_v \theta_{v-1}^p \sin(p\theta_v),$$

and it is easy to prove that

$$J_{u,v} = -2\mathcal{S}_u\mathcal{S}_v \frac{\theta_{v-1}\sin\theta_v}{1 - 2\theta_u\theta_{v-1}\cos\theta_v + \theta_u^2\theta_{v-1}^2}.$$

When $u,v \in \mathcal{P}_\mu \bigcup \mathcal{Z}_\mu \bigcup \mathcal{P}_\phi \bigcup \mathcal{Z}_\phi$, we can apply the formulas given in [5] for the computation of $J_{u,v}$ in case all the poles and the zeros are purely complex.

Analyzing the sign of the product $\mathcal{S}_u\mathcal{S}_v$, we find that the matrix $\mathbf{J}(\theta)$ can be re-written more compactly as $\mathbf{J}(\theta) = \begin{bmatrix} \mathbf{G} & -\mathbf{C} \\ -\mathbf{C}^\top & \mathbf{H} \end{bmatrix}$, where the size of the block matrix $\mathbf{C}$ is $n \times m$. The identity $\left| \begin{smallmatrix} \mathbf{G} & -\mathbf{C} \\ -\mathbf{C}^\top & \mathbf{H} \end{smallmatrix} \right| = \left| \begin{smallmatrix} \mathbf{G} & \mathbf{C} \\ \mathbf{C}^\top & \mathbf{H} \end{smallmatrix} \right|$ leads to the conclusion that $\int_\Theta |\mathbf{J}(\theta)|^{1/2} d\theta$ has the same value for the models ARMA(n,m), ARMA(n+m,0), ARMA(0,n+m). A similar conclusion was drawn in [4] for the particular case when all the poles and the zeros are real-valued.

*Estimation of AR and ARMA models by stochastic complexity* 59[15] RISSANEN, J. (1978). Modeling by shortest data description. *Automatica* **14** 465–471.
[16] RISSANEN, J. (1986). Order estimation by acumulated prediction errors. *J. Appl. Prob.* **23A** 55–61. MR0803162
[17] RISSANEN, J. (1989). *Stochastic Complexity in Statistical Inquiry*. World Scientific Publ. Co., River Edge, NJ, 175 pp. MR1082556
[18] RISSANEN, J. (1996). Fisher information and stochastic complexity. *IEEE Trans. Inf. Theory* **42** (1, Jan.) 40–47. MR1375327
[19] RISSANEN, J. (2000). MDL denoising. *IEEE Trans. Inf. Theory* **46** (7, Nov.) 2537–2543.
[20] SCHWARZ, G. (1978). Estimating the dimension of the model. *Ann. Stat.* **6** 461–464. MR0468014
[21] SEGHOUANE, A.-K. AND BEKARA, M. (2004). A small sample model selection criterion based on Kullback's symmetric divergence. *IEEE Trans. Signal. Proces.* **52** 3314–3323. MR2107913
[22] WAX, M. (1988). Order selection for AR models by predictive least squares. *IEEE Trans. on Acoustics, Speech and Signal Processing* **36** 581–588.